\documentclass[12pt]{article}  %% generate by rtf2latex

    \usepackage{graphicx}
    \usepackage{amsmath}
    \usepackage{amsthm}
    \usepackage{amssymb}
    \usepackage{latexsym}

\theoremstyle{definition}

\newcommand{\cF}{\mbox{$\mathcal{F}$}}
\newcommand{\BR}{\mbox{$\mathbb{R}$}}
\newcommand{\BN}{\mbox{$\mathbb{N}$}}

\title{\bf Bounded linear operators in PN-spaces}
\author{ Delavar Varasteh Tafti$^1$, Mahdi Azhini$^2$}
\date{
Department of Mathematics, Science and Research Branch,
 Islamic Azad University, Tehran, Iran.}
\begin{document}
\maketitle
\begin{abstract}
The manner in  which the  strong completeness of probabilistic Banach spaces
is frequently exploited  depends on the Baire theorem about strong complete
probabilistic metric spaces. By using the Baire theorem, in this paper,
 we present the open mapping, closed Graph, Principle of Uniform Boundedness
and Banach-Steinhaus theorems in PN-spaces.
%%%%

\noindent
{\bf Key words:}
$\check{S}$erstnev probabilistic normed spaces, bounded linear operators, the
open mapping theorem.\\
{\bf Mathematics Subject Classification (2010):} 46S50, 54E70.
\footnote{\noindent $^1$
 Email: delavar.Varasteh@gmail.com}
\footnote{ $^2$ Corresponding author:
Mahdi Azhini. Email:m.azhini@srbiau.ac.ir}
\end{abstract}
%%%%%
%%%%%%%
\section*{1. Introduction}
Probabilistic Normed (briefly PN) spaces  were first  introduced by
$\check{S}$erstnev in  a series of papers [1-4]. Then a new definition was
proposed by Alsina, Schweizer and Sklar [2]. Linear operators in probabilistic
normed spaces were first studied by Lafuerza-Guill$\acute{e}$n,
Rodrignez-Lallena, and Sempi in [5]. In section 2, we recall some notations
and definitions of probabilistic normed spaces according to those of  [2] and
[8].  In section 3, we present  another way of formulating a PN-space. In
section 4, by using this formulating, we prove all the known classical
theorems of functional analysis, such as Open Mapping and Closed Graph and
Principle of Uniform Boundedness and Banach-Steinhaus theorems in PN-spaces.
%%%%%%%%%%%%%%
\section*{2. Preliminaries}
%%%%%%%%%
\paragraph{Definition 2.1.}
A  distribution function (briefly a d.f.) is a function
$F:\bar{\BR}\longrightarrow  [0,1]$
that is nondecreasing and  left-continuous on $\BR$, Moreover, $F(-\infty)=0$
and $F(+\infty)=1$. Here $\bar{\BR}=\BR\cup\{-\infty,+\infty\}$. The set of
all the d.f.'s will be denoted by $\Delta$ and  the subset of those d.f.'s
called distance d.f.'s, such that $F(0)=0$, by $\Delta^+$. We shall also
consider $\mathcal{D}$ and $\mathcal{D}^+$,
 the subsets of $\Delta$ and $\Delta^+$,
respectively, formed by the proper d.f.'s, i.e., by those d.f.'s $F\in\Delta$
that  satisfy the conditions
$$\lim_{x\longrightarrow  -\infty}F(x)=0\quad \text{and} \quad
\lim_{x\longrightarrow
+\infty}F(x)=1.$$
The first of these  is obviously satisfied  in all of $\Delta^+$,
 since, in it,
$F(0)=0$. By setting $F\leq G$ whenever $F(x)\leq G(x)$ for every $x\in\BR$,
one introduces in natural ordering in $\Delta$ and in $\Delta^+$. The maximal
element for $\Delta^+$ in this order in the d.f., given by
$$H_a(x)=\begin{cases}
0 & \text{if }\quad  x\leq a,\\
1 & \text{if} \quad x>a
\end{cases}$$
where $a\in\BR$.

The space $\Delta$
 can be metrized in several way [10], [12], [14], but we shall here
adopt the modified Levy metric $d_L$. If $F$ and $G$ are d.f.'s and $h$ is in
$(0,1]$, let $(F, G;h)$ denote the condition
$$F(x-h)-h\leq G(x)\leq F(x+h)+h\quad \text{for all}~ x\in
(-\frac{1}{h},\frac{1}{h}).$$
Then the modified Levy metric $d_L$ is defined by
$$d_L(F,G)=\inf\{h\in (0,1]:\text{both}(F,G;h)\text{and}(G,F;h)\text{hold}\},$$
under which a sequence of distribution functions $\{F_n\}_{n\in \BN}$
converges to $F\in\Delta^+$ if and only if at each point $t\in\BR$, where $F$
is
continuous, $F_n(t)\longrightarrow  F(t)$ [7].
%%%%
\paragraph{Definition 2.2.}
A triangle function is a  binary operation on $\Delta^+$, namely a function
$\tau:\Delta^+\times \Delta^+\longrightarrow  \Delta^+$ that is associative,
commutative,
nondecreasing and which has $H_0$ as unit, that is, for all $F,G,H\in\Delta^+$,
we have
\\
$\tau(\tau(F,G),H)=\tau(F,\tau(G,H)),$\\
$\tau(F,G)=\tau(G,F),$\\
$\tau(F,H)\leq \tau(G,H), \quad \text{if}\; F\leq G,$\\
$\tau(F,H_0)=F.$
\\
Typical continuous triangle function are convolution and  the operators
$\tau_T$ and $\tau_{T^*}$, which are,  respectively, given by
\\
$ \tau_T(F,G)(x)= \sup\{T(F(s),G(t))|s+t=x\}\quad \text{and}$\\
$ \tau_{T^*}(F,G)(x)= \inf\{T^*(F(s),G(t))|s+t=x\}$
\\
for all $F,G$ in $\Delta^+$ and all $x\in \BR$ [7], here $T$ is a continuous
$t$-norm, i.e, a continuous binary operator on $[0,1]$ that is associative,
commutative, nondecreasing and has 1 as identity; $T^*$ is a continuous
$t$-conorm, namely a continuous binary operation on $[0,1]$ that is related to
continuous $t$-norm through
$$T^*(x,y)=1-T(1-x, 1-y).$$
The most important $t$-norms  are the functions $W,\Pi$ and $M$ which are
defined, respectively by
$$\begin{array}{l}
W(a,b)=\max\{a+b-1,0\},\\
\Pi(a,b)=a.b,\\
M(a,b)=\min\{a,b\}.
\end{array}$$
Their corresponding $t$-conorms are given, respectively by
$$\begin{array}{l}
W^*(a,b)=\min\{a+b,1\},\\
\Pi^*(a,b)=a+b-ab,\\
M^*(a,b)=\max\{a,b\}.
\end{array}$$
It follows without difficulty that
$$\tau_T(H_a,H_b)=H_{a+b}=\tau_{T^*}(H_a,H_b),$$
for any continuous $t$-norm $T$, any continuous $t$-conorm $T^*$ and any
$a,b\geq 0$.
%%%%%
\paragraph{Definition 2.3.} [5]
A $\check{S}$erstnev PN-space or a $\check{S}$erstnev space is a triple
$(V,\nu,\tau)$, where $V$ is a (real or complex) vector space, $\nu$ is a
mapping from $V$ into $\Delta^+$ and $\tau$ is a continuous triangle function
and
the following conditions are  satisfied for all $p$ and $q$  in $V$:
\\
$(N1)\
 \nu_p=H_0, \; \text{if and only if,}~  p=\theta (\theta~ \text{is the
null vector in} V);$\\
$(N3) \ \nu_{p+q}\geq \tau(\nu_p,\nu_q);$ \\
$(\check{S})$
$\  \forall  \alpha\in\BR\backslash\{0\}, \forall x\in \bar{\BR}_+,
\nu_{\alpha p}(x)=\nu_p(\frac{x}{|\alpha|}).$\\
$\  \text{Notice that condition} (\check{S}) \text{implies}$\\
$(N2) \ \forall p\in V, \nu_{-p}=\nu_p.$
%%%%%%
\paragraph{Definition 2.4.} [5]
A PN-space is a quadruple $(V,\nu,\tau,\tau^*)$, where $V$ is a real vector
space, $\tau$ and $\tau^*$ are continuous triangle functions such that
$\tau\leq \tau^*$, and the mapping $\nu:V\longrightarrow  \Delta^+$ satisfies,
for all $p$
and $q$ in $V$, the conditions:
\\
$(N1)\
\nu_p=H_0, \; \text{if and only if,} p=\theta (\theta\text{ is the null
vector in V)};$\\
$(N2)\  \forall p\in V, \; \nu_p=\nu_p;$\\
$(N3)\ \nu_{p+q}\geq \tau(\nu_p,\nu_q);$\\
$(N4)\  \forall \alpha\in [0,1], \nu_p\leq \tau^*(\nu_{\alpha p},
\nu_{(1-\alpha)p}).$
\\
The function $\nu$ is called the probabilistic norm.

If $\nu$ satisfies the condition, weaker than $(N1), \nu_\theta=H_0$, then
$(V,\nu,\tau,\tau^*)$ is called a Probabilistic Pseudo-Normed space (briefly,
a PPN space).

If $\nu$ satisfies the condition $(N1)$ and $(N2)$, then $(V,\nu,\tau,\tau^*)$
is said  to be a Probabilistic seminormed space (briefly, PSN space).

If $\tau=\tau_T$ and $\tau^*=\tau_{T^*}$ for some continuous $t$-norm $T$ and
its $t$-conorm $T^*$, then $(V,\nu,\tau_T,\tau_{T^*})$ is denoted by
$(V,\nu,T)$ and is called a Menger space.
%%%%\
\paragraph{Definition 2.5.}[5]
There is a natural topology on a PN-space $(V,\nu,\tau,\tau^*)$, called the
strong topology, which is defined by the system of neighborhoods
$$N_p(t)=\{q\in V|~ \nu_{p-q}(t)>1-t\}=\{q\in
V|d_L(\nu_{p-q},H_0)<t\},\text{where}~ p\in V\text{and}~ t>0.$$
Let $(V,\nu,\tau,\tau^*)$ be a PN-space, then
\begin{itemize}
\item[i)]
A sequence $\{p_n\}$ in $V$ is said to be strongly convergent to a point $p$
in $V$, and we write $p_n\longrightarrow  p$, if for each $t>0$, there exist  a
positive
integer $m$ such that  $p_n\in N_p(t)$, for $n\geq m$.
\item[ii)]
A sequence $\{p_n\}$ in $V$ is called a strong cauchy sequence  if for every
$t>0$, there is a positive integer $N$ such that $\nu_{p_n-p_m}(t)>1-t$,
whenever  $m,n>N$.
\item[iii)]
The PN-space $(V,\nu,\tau,\tau^*)$ is said to be distributionally compact
($D$-compact) if every sequence $\{p_n\}$ in $V$  has a convergent subsequence
$\{p_{n_k}\}$. A subset $A$ of  a PN-space $(V,\nu,\tau,\tau^*)$ is said to be
$D$-compact if every sequence $\{p_m\}$ in $A$ has a subsequence $\{p_{m_k}\}$
that converges to a point  $p\in A$.
\item[iv)]
In the strong topology, the closure  $\overline{N_p(t)}$ of $N_p(t)$ is
defined  by $\overline{N_p(t)}=N_p(t)\cup N'_p(t)$, where $N'_p(t)$ is the set
of limit  points of all convergent sequence in $N_p(t)$.
\end{itemize}
In this paper, we will consider  those PN-space for which
$\tau=\tau_M=\tau^*$, which are clearly $\check{S}$erstnev as well as Menger
PN-space also, and will  denote such spaces simply by the pair $(V,\nu)$.
%%%%%%%%%%%%%%%%%%%
\section*{3. Bounded linear operators in SPN-spaces}
\paragraph{Definition 3.1.}
Let $\nabla$ be the set of all nondecreasing and right continuous functions
$f:[0,1]\longrightarrow  \bar{\BR}$, and $\nabla^+$ the  subset consisting of
all
non-negative $f\in \nabla$. Suppose $\mathcal{R}^+$in the set of all $f\in
\nabla^+$ with $f(w)<+\infty$, for all $w\in (0,1)$. For $F\in \Delta$, define
$\hat{F}:[0,1]\longrightarrow  \BR$ by
$$\hat{F}(w)=\sup\{t\in \BR| F(t)<w\}.$$
%%%%%
\paragraph{Theorem 3.2.}[7]
If $F\in \mathcal{D}^+$, then $\hat{F}\in \mathcal{R}^+$. Moreover, for $F$ and
$G$
in
$\Delta^+$,
\begin{itemize}
\item[i)]
if $F\leq G$, then $\hat{F}\geq \hat{G}$;
\item[ii)]
if $\hat{F}\geq \hat{G}$, then  $F(t)\leq G(t+h)$ for all  $t\in \BR$ and
$h>0$;
\item[iii)]
$\widehat{(\tau_M(F,G))}=\hat{F}+\hat{G}$;
\item[iv)]
$\widehat{F(\frac{t}{h}))}=h\hat{F}(t), \text{for all } t\in \BR~ \text{and}
~h>0.$
\end{itemize}
%%%%%%%
\paragraph{Theorem 3.3.}[7]
The map $\wedge:\Delta^+\longrightarrow  \nabla^+$ is one-to-one.
%%%%%
\paragraph{Definition 3.4.}
Let $(V,\nu)$ be a PN-space. The map $V\longrightarrow  \mathcal{R}^+$ given
through
the
composition of maps $V\overset{\nu}{\longrightarrow }
\Delta^+\overset{\widehat{ }}{\longrightarrow}\mathcal{R}^+$
will be  denoted by $\|\cdot\|$, and for $x\in V$, the value of $\|x\|$ at
$w\in (0,1)$, simply by $\|x\|_w$, i.e, $\|x\|_w=(\widehat{\nu_x})(w)$.
The proof of the following Theorem is obtained, using the definition of
PN-space, Theorem 3.2.
%%%%%
\paragraph{Theorem 3.5.}
Suppose $(V,\nu)$ is a PN-space. The map $\|\cdot\|:V\longrightarrow  \nabla^+$
has the
following properties:
\begin{itemize}
\item[i)]
For $x\in V$, $\|x\|\geq 0$, and $\|x\|=0$ if and only if $x=0$.
\item[ii)]
For $x\in V$ and $\alpha\in\BR$,  $\|\alpha x\|=|\alpha|\|x\|.$
\item[iii)]
For $x$ and $y$ in $V$, $\|x+y\|\leq \|x\|+\|y\|.$
\end{itemize}
%%%%%%%
\paragraph{Remark 3.6.}
Note that all assertion of theorem 3.5 are to be understood point-wise. Thus
for every $w\in (0,1)$, the map $\|\cdot\|_w:V\longrightarrow  \BR$ is a norm
on $V$.
%%%%%%
\paragraph{Theorem 3.7.}[12]
Let $(V,\nu)$ be a PN-space. For $p\in V, r>0$, and $w\in (0,1)$ let
$B_w(p;r)$ be defined by
$B_w(p;r)=\{x\in V| \|x-p\|_w<r\}$.  Then the family
$\{B_w(p;r)|p\in V, r>0, w\in (0,1)\}$ forms a basis for the strong topology
on $V$.
%%%
\paragraph{Theorem 3.8.}[6]
A PN-space $(V,\nu)$ is a topological vector space if and only if  $\|x\|\in
\mathcal{R}^+$ for every $x\in V$. Moreover it is seen this topology is locally
convex.\\[4mm]
Now, we denote  the strong topology of PN-space $(V,\nu)$ by $\tau$, and  the
corresponding dual space by $V^*$. Hence $f\in V^*$  is a continuous
(equivalently bounded because of PN-space $(V,\nu)$ (see [15])) linear
functional on $V$.
%%%%%%
\paragraph{Theorem 3.9.}[12]
 Let $f:V\longrightarrow \BR$ be linear. Then $f\in V^*$ if and only if there
is some $w\in
(0,1)$, with
$$\sup_{x\in B_{V,1-w}}|f(x)|<+\infty, \text{where } B_{V,1-w}=\{x\in
V|\|x\|_w\leq 1\}.$$
Moreover if $\|f\|:[0,1)\longrightarrow  [0,+\infty)$ is defined by
$\|f\|_w=\inf_{w<w'}(\sup_{x\in B}|f(x)|)$, then $\|f\|$ is a nondecreasing
and right-continuous function, and
$$|f(x)|\leq \|f\|_w\|x\|_{w'}, \text{for all}~
 x\in V \text{and}~ w\in (0,1).$$
By theorem 3.9, a linear functional $f:V\longrightarrow  \BR$ belongs to $V^*$
if and only
if $\|f\|_w<+\infty$ for some $w\in (0,1)$.
%%%%%%%
\\[4mm]
The following theorem introduces a simple way of constructing  PN-space. We
recall that for a nondecreasing function $f:(0,1)\longrightarrow  \BR$, the
function
$l^-f:(0,1)\longrightarrow  \BR$ is defined by
$$l^-f(w_0)=\lim_{w\longrightarrow  w_0^-}f(w)=\sup_{w<w_0}f(w).$$
Clearly $l^-f$ is nondecreasing and left-continuous. We also denote the
Lebesgue measure on $(0,1)$ by $m$.
%%%%%%%%%%%%%%%
\paragraph{Theorem 3.10 [13]}
Let $V$ be a real vector space. Suppose  $p:V\times (0,1)\longrightarrow
[0,+\infty)$,
satisfies the following conditions:
\begin{itemize}
\item[i)]
for every $x\in V$, the function $P(x,.):(0,1)\longrightarrow  [0,+\infty)$ is
nondecreasing,
\item[ii)]
for every $w\in (0,1)$, the map $p(.,w):V\longrightarrow  [0,+\infty)$ is a
semi-norm on
$V$,
\item[iii)]
the family of semi-norms $\{p(.,w)|w\in (0,1)\}$ is separating on $V$.
\end{itemize}
then, there exists a unique map $\nu:V\longrightarrow \Delta^+$, defined by
$$\nu_x(t)=m(\{w\in (0,1)|p(x,w)<t\}),$$
such that $(V,\nu)$ is a PN-space. Moreover, $\|x\|_w=l^-p(x,w)$ for all $x\in
V$ and $w\in (0,1)$.
%%%%%
\paragraph{Theorem 3.11. [13] }
Let $V$ be a real vector space. Suppose $p:V\times(0,1)\longrightarrow
  [0,+\infty)$
is a map
which satisfies conditions of theorem 3.10, with $\nu:V\longrightarrow
\Delta^+$
the
corresponding probabilistic norm on $V$. Then the strong topology on $V$ is
induced by the  separating family of semi-norm
$\{p (.,w)|w\in (0,1)\cap Q\}$. hence $V$ under the strong topology is a
locally topological vector space.
%%%%
\paragraph{Remark 3.12.}
Let $(V,\nu)$ be a PN-space, and $p:V\times (0,1)\longrightarrow  [0,+\infty)$
be given by
$p(x,w)=\|x\|_w=\hat{\nu}_x(t)$. Then by Remak 3.6 and Theorem 3.11, the
family
$$\{B_{V,w}(0;r)|w\in (0,1), r>0\},$$
where $B_{V,w}(0;r)$ is given by $\{x\in V, \|x-x_0\|_w<r\}$, forms a local
basis for the strong topology. Moreover, a subset $E$ of $V$ is bounded if and
only if the map $p (.,w):V\longrightarrow  [0,+\infty)$ is bounded on $E$, for
every $w\in
(0,1)$ (see [16] ).

Using theorem 3.10, we may construct many examples of PN-spaces.
%%%%%%%
\paragraph{Example 3.13.}
Let $(V,\|\cdot\|)$ be a real normed vector space. If we define
$p:V\times(0,1)\longrightarrow  [0,+\infty)$ by
$$p(x,w)=\|x\|\quad \forall x\in V, \quad \forall w\in (0,1),$$
then $p$ satisfies all condition of theorem 3.10. Using $\nu_x(t)=m(\{w\in
(0,1)|p(x,w)<t\})$, we obtain $\nu_x=H_{\|x\|}$.
%%%%
\paragraph{Example 3.14.}
Let $(V,\nu)$ and $(W,\mu)$ be two PN-spaces. Then, the map $\theta:V\times
W\longrightarrow
\Delta^+$ defined by
$$\theta(x,y)=\tau_M(\nu_x,\nu_y),$$
is a probabilistic norm on $V\times W$, i.e., $(V\times W,\theta)$ is a
PN-space. Because define $p:V\times W\times(0,1)\longrightarrow  [0,+\infty)$
by
$$p(x,y,w)=\|x\|_{V,w}+\|y\|_{W,w}$$
where $\|\cdot\|_V:V\longrightarrow  \mathcal{R}^+$
 and $\|\cdot\|_W:W\longrightarrow
\mathcal{R}^+$ are
the composite
maps given by Definition 3.4. Then the desired is obtained by Theorem 3.10 and
Theorem 3.2 (iii).
%%%%%%%%%%%%%%%%5
\paragraph{Definition 3.15.}
Let $V$ and $W$ be two PN-spaces. We remind that, for topological vector space
$V$ and $W$, a linear operator $T:V\longrightarrow  W$ is called bounded if
$T(E)$ is
bounded (in $W$) whenever $E$ is bounded subset of $V$. Metrisability of the
PN-spaces implies that $T$ is bounded if and only if $T$ is continuous (see
[16]
). As usual, the set of all bounded linear operators from $V$ to $W$ is
denoted by $B(V,W)$.
%%%%%%
\paragraph{Definition 3.16.}
Let $V$ and $W$ be two PN-spaces.
For linear operator $T:V\longrightarrow  W$ and  a pair
$(w,w')\in J=\{(w,w')|w,w'\in(0,1)\}\subseteq \BR^2$,
the norm $\|T\|_{(w,w')}$ is defined by
$$\|T\|_{(w,w')}=\sup\{\|Tx\|_{w'}|~x\in B_{V,w}\}, \text{where}
~B_{V,w}=\{x\in V| \|x\|_w\leq 1\}.$$
%%%%%
\paragraph{Theorem 3.17.} [13]
For a pair $(w,w')\in J$, if $\|T\|_{(w,w')}<+\infty$, then
$$\|Tx\|_{w'}\leq \|T\|_{(w,w')}\|x\|_{w}, \text{for all}~ x\in V.$$
%%%%%%
\paragraph{Theorem 3.18.} [13]
Let $V$ and $W$ be two PN-spaces, and $T:V\longrightarrow  W$ be a linear
operator. Then
$T\in B(V,W)$ if and only if, for every $w'\in (0,1)$ there exists $w\in
(0,1)$ such that $\|T\|_{(w,w')}<+\infty$.
%%%%%
%%%%%%%%%%%%%%%%%%%%%
\section*{4. Main Results}
Baire theorem is one of the two keystones of Linear Functional Analysis, the
other one  being the Hahn-Banach theorem. Bair's theorem for reaching
consequences include such basic theorems as the Banach open mapping theorem,
the Banach closed graph theorem and the Banach-Steinhaus theorem. This chapter
concludes consequences of Bair's theorem in PN-spaces.

It is impossible to speak about PN-spaces without making reference to concept
of a Probabilistic Metric space (briefly a PM space).
PM spaces are studied in depth in the fundamental book [7] by Schweizer and
Sklar.
%%%%%%%
\paragraph{Definition 4.1.}
A triple $(V,\cF,\tau)$, where $V$ is a nonempty set and $\cF$ is a function from
$V\times V$ into $\Delta^+$ and $\tau$  is a triangle function, is a
Probabilistic Metric space if the following  conditions are satisfied for all
points $p,q$ and $r$ in $V$:\\
$(PM1) \ {\cF}(p,p)=H_0;$\\
$(PM2)\ {\cF}(p,q)\neq H_0~ \text{if}~ p\neq q;$\\
$(PM3)\ {\cF}(p,q)={\cF}(q,p);$\\
$(PM4)\ {\cF}(p,r)\geq \tau({\cF}(p,q),{\cF}(q,r)).$
\\
Notice that, given PN space $(V,\nu)$, a probabilistic metric ${\cF}:V\times
V\longrightarrow  \Delta^+$ is defined via ${\cF}(p,q)=\nu_{p-q}$, so that the
triple
$(V,{\cF},\tau)$ becomes a PM space.
%%%%%
\paragraph{Theorem 4.2 (Baire Theorem in PM spaces).[11]}
Let $(V,{\cF},\tau)$ be a strong complete PM space under a continuous triangle
function $\tau$. If $\{G_n\}$ is a sequence of dense and strongly open subset
of $V$, then $\cap_{n=1}^\infty G_n$ is not empty. (In fact it is dense in
$V$).
%%%%%
\paragraph{Corollary 4.3.}
Let $(V, {\cF},\tau)$ be a strong complete PM space under a continuous
triangle function $\tau$. Then the  following two equivalent properties hold:
\begin{itemize}
\item[a)]
Let $\{f_n\}_{n=0}^\infty$  be a sequence of strongly closed subsets of $V$
such that $\text{int}
 F_n=\phi$ for all $n\geq 0$. Then $\text{int}(\cup_{n=0}^\infty
F_n)=\phi$.
\item[b)]
Let $\{G_n\}_{n=0}^\infty$ be a sequence  of strongly open subsets of $V$
such that $G_n=V$ for all $n\geq 0$. Then $\overline{\cap_{n=0}^\infty
G_n}=V$.
\end{itemize}
%%%%
\paragraph{Corollary 4.4.}
Let $(V,{\cF},\tau)$ be a PM space under a continuous triangle function
$\tau$, and let  $F_n$, $n\geq 0$, be strongly closed subsets of $V$ such that
 $V=\cup_{n=0}^\infty F_n$.
\begin{itemize}
\item[a)]
 If $\text{int} F_n=\phi$ for all $n\geq 0$, then $V$ is not strong complete.
\item[b)]
If $V$ is strong complete, there exists $n_0\geq 0$  such that  $\text{int}
F_{n_0}\neq \phi$.
\end{itemize}
%%%%
\paragraph{Theorem 4.5 (Open Mapping theorem in PN-spaces).}
If $T$ is a probabilistic bounded linear operator from a strong complete
PN-space $(V,\nu)$ onto a strong complete $(W,\mu)$, then $T$ is an open
mapping.
%%%%
\paragraph{Proof.}
The Theorem will be proved by the following steps: \\
{\bf Step 1:} The set $T(B_{V,w}(0;1))$, where $B_{V,w}(0;r)=\{x\in
V|\|x\|_w<r\}$ for all  $r>0$ and $w\in (0,1)$, contains an strong open ball.
Now fix $w\in (0,1)$. Given any $y\in W$,
%%%%%%%%%%%%%%%%%%pag 15
%%%% page 15
there exists $x\in V$ such that $y=Tx$, since $T$ is sujective. Since $x\in
B_{V,w}(0;n)$ for some integer $n\geq 1$, this shows that
$W=\bigcup_{n=1}^\infty T(B_{V,w}(0;n))$, hence
$W=\bigcup_{n=1}^\infty \overline{T(B_{V,w}(0;n))}$.
The space $W$  being strong complete, hence Corollary 4.4 (b) shows that there
exists  an integer $n_0\geq 1$ such that
$$\text{int}\overline{T(B_{V,w}(0;n_0))}\neq \phi.$$
Therefore $\text{int}\overline{T(B_{V,w}(0;1))}\neq \phi$, since
$\overline{T(B_{V,w}(0;1))}=\frac{1}{n_0}\overline{T(B_{V,w}(0;n_0))}$
by the linearity of $T$. Hence the set $\overline{T(B_{V,w}(0;1))}$ contains
on strong open ball. \\
{\bf Step 2:}
The set $\overline{T(B_{V,w}(0;1))}$ contains an strong open ball centered at
the origin of $W$. By Step 1, there exist $y\in W$ and $s>0$ such that
$B_{W,w}(y;2s)\subset \overline{T(B_{V,w}(0;1))},$ and hence
$$B_{W,w}(0;2s)=\{-y\}+B_{W,w}(y;2s)\subset
\{-y\}+\overline{T(B_{V,w}(0;1)))}.$$
Since $-y\in \overline{T(B_{V,w}(0;1))}$(because $y\in
\overline{T(B_{V,w}(0;1))}$ and $T$ is linear), it follows that
$$\{-y\}+\overline{T(B_{V,w}(0;1))}\subset 2T(B_{V,w}(0;1).$$
{\bf Step 3:}
The set $T(B_{V,w}(0;1))$ contains an strong open ball centered at the origin
of $W$. To prove this assertion, we will show that
$B_{W,w}(0;\frac{s}{2})\subset T(B_{V,w}(0;1))$, where $s>0$ is the radius of
the ball $B_{W,w}(0;s)$ found in Step 2. This means that, given any $y\in
B_{V,w}(0;\frac{s}{2})$, we need to find $x\in B_{V,w}(0;1)$ such that $y=Tx$.
So, let $y\in B_{W,w}(0;\frac{s}{2})$ be given. Since $y\in
B_{W,w}(0;\frac{s}{2})\subset \overline{T(B_{V,w}(0;\frac{1}{2}))}$ by Step 2,
there exists $x_1\in B_{V,w}(0;\frac{1}{2})$ such that
$\|y-Tx_1\|_w<\frac{s}{2^2}$, since
$y-Tx_1\in B_{W,w}(0;\frac{s}{2^2})\subset
\overline{T(B_{V,w}(0;\frac{1}{2^2}))},$ again by Step 2, thee exists
$x_2\in B_{V,w}(0;\frac{1}{2^2})$ such that $\|y-Tx_1-Tx_2\|_w<\frac{s}{2^3}$,
 and so on. In this manner we construct a sequence $\{x_n\}_{n=1}^\infty$ of
points $x_n\in V$ with the following properties:
$x_n\in B_{V,w}(0;\frac{1}{2^n})$ and
$\|y-T(\sum_{k=1}^nx_k)\|_w<\frac{s}{2^{n+1}}$, for all $n\geq 1$
.
Since the series $\sum_{n=1}^\infty x_n$ is absolutely convergent (because
$\|x_n\|_w<\frac{1}{2^n}$, for all $n\geq 1$) and the space $V$ is strong
complete, the series $\sum_{n=1}^\infty x_n$ converges to a point $x\in V$ and
$$\|x\|_w\leq \sum_{n=1}^\infty \|x_n\|_w<\frac{1}{2}+\frac{1}{2^2}+\cdots +
\frac{1}{2^{n+1}}+\cdots=1,$$
by Theorem 3.5, so that $x\in B_{V,w}(0;1)$. Furthermore,
$y=\lim_{n\longrightarrow \infty}T(\sum_{k=1}^n x_k)=Tx,$ since  $T$ is
continuous.  Hence
the assertion is proved. \\
{\bf Step 4:}
The mapping $T$ is open strongly.
 Given any strong open subset $U$ of $V$ and give any
$y\in T(U)$, we must find $\delta>0$ such that $B_{W,w}(y;\delta)\subseteq
T(U)$.
So let $x\in U$ be such that $y=Tx$. Since $U$ is strong  open, thee exists
$r>0$ such that $B_{V,w}(x;r)=\{x\}+B_{V,w}(0;r)\subseteq U$, and by Step 3,
there exists $\delta>0$ such that $B_{W,w}(0;\delta)\subseteq T(B_{V,w}(0;r))$.
Hence
\begin{align*}
B_{W,w}(y;\delta)= \{y\}+B_{W,w}(0;\delta) & \subseteq \{y\}+T(B_{V,w}(0;r))\\
& = T(\{x\}+B_{V,w}(0;r))\\
& \subseteq T(U).
\qquad \qquad \qquad
\qquad \qquad \qquad
\qquad \qquad \qquad
 \Box
\end{align*}
%%%%%%
\paragraph{Corollary 4.6.}
Let $(V,u)$ and $(W,\mu)$ be two  strong complete PN-spaces and let $T\in
B(V,W)$ be bijective. Then $T^{-1}\in B(W,V)$.
%%%%%%
\paragraph{Proof.}
First fix $w\in (0;1)$. It is clear that $T^{-1}:W\longrightarrow  V$ is
linear, since
$T\in B(V,W)$, hence by Definition 3.15, suffices to show that
$T^{-1}:W\longrightarrow  V$
is continuous at the origin, i.e, that given any strong open ball
$B_{V,w}(0;\epsilon)\subseteq V$, there exists an strong  open ball
$B_{V,w}(0;\delta)\subseteq W$, such that $T^{-1}(B_{W,w}(0;\delta))\subseteq
B_{V
,w}(0;\delta)$, by definition of continuouty at a point. But the inclusion
$T^{-1}(B_{W,w}(0;\delta))\subseteq B_{V,w}(0;\epsilon)$, is equivalent to the
inclusion $B_{W,w}(0;\delta)\subseteq T(B_{V,w}(0;\epsilon))$, since the
mapping
 $T:V\longrightarrow  W$ is bijective, and open Mapping Theorem precisely show
that this
last inclusion hold for some $\delta>0$.
\hfill $\Box$
%%%%%%
\paragraph{Theorem 4.7.}
Let $\|\cdot \|_w$ and $\|\cdot\|'_w$, for all $w\in (0,1)$,  be two norms on
the same vector space $V$, with the following properties: both
$(V,\|\cdot\|_w)$ and $(V,\|\cdot\|'_w)$, are strong complete
PN-spaces, and there exists
a constant $C$ such that  $\|x\|'_w\leq C\|x\|_w$ for all $x\in V$. Then there
exists  $w'\in (0,1)$ and $c_{w'}\in\BR$ such that
$\|x\|_{w'}\leq c_{w'}\|x\|'_w$, for all $x\in V$.
%%%%%%
\paragraph{Proof.}
Fix $w\in (0,1)$. The bijective and linear identity mapping
$i:(V,\|\cdot\|_w)\longrightarrow  (V,\|\cdot\|'_w)$ belongs to $B(V,V)$ by
assumption.
Corollary 4.6, therefore shows that  the inverse mapping
$i^{-1}:(V,\|\cdot\|'_w)\longrightarrow  (V,\|\cdot\|_w)$, belongs to $B(V,V)$,
hence by Theorems 3.18 and 3.17,
there exists $w'\in (0,1)$ and $c_{w'}$ such that  $\|x\|_{w'}\leq
c_{w'}\|x\|'_w$, for all $x\in V$.
\hfill $\Box$
%%%%%%
\paragraph{Definition 4.8.}
Let $(V,\nu)$ and $(W,\mu)$ be two PN-spaces. The graph $GrT$ of a mapping
$T:V\longrightarrow  W$ is the subset of the product $V\times W$ defined by
$$GrT=\{(x,Tx)\in V\times W; x\in V\}.$$
Let $(V,\nu)$ and $(W,\mu)$ be two PN-spaces. By Example 3.14 $(V\times
W,\theta)$ is a PN-space, where $\theta(x,y)=\tau_M(\nu_x,\mu_y)$, for all
$x\in V$ and $y\in W$, is a probabilistic norm on $V\times W$.
%%%%%
\paragraph{Definition 4.9.}
If $(V,\nu)$ and $(W,\mu)$ are PN-spaces, a mapping $T:V\longrightarrow  W$ is
said to be
closed graph, if its graph $GrT$ is strong closed in the PN-space $V\times W$
(equipped with strong product topology, by Example 3.14 and  Theorem 3.11).
Therefore a mapping $T:V\longrightarrow  W$ is strongly closed if and only if
$x_n\longrightarrow  x$
in $V$ and
$Tx_n\longrightarrow  y$ in $W$ implies $y=Tx$.
%%%%%
\paragraph{Theorem 4.10 (Closed Graph Theorem in PN-spaces).}
Let $(V,\nu)$ and $(W,\mu)$ be strongly complete PN-spaces, and
$T:V\longrightarrow  W$ be a
closed linear operator. Then $T\in B(V,W)$.
\paragraph{Proof.}
 First fix $w\in (0,1)$ and
$w'\in(0,1)$. By using Theorem 3.5, we define another norm on $V$ by
$$\|x\|'_{V,w}=\|x\|_{V,w}+\|Tx\|_{W,w'},$$
and $x\in V$. On the other hand by the definition of norm $\|\cdot\|'_{V,w}$,
if $\{x_n\}_{n=1}^\infty$ is a Cauchy sequence with respect to
$\|\cdot\|'_{V,w}$, then $\{x_n\}_{n=1}^\infty$ and $\{Tx_n\}_{n=1}^\infty$
are Cauchy sequences in the PN-spaces $V$ and $W$, respectively.
Since both spaces are strongly
 complete, there exist $x\in V$ and $y\in W$ such
that $x_n\longrightarrow  x$ in $V$ and $Tx_n\longrightarrow  y$ in $W$, and
thus $y=Tx$,
since $T$ is strongly
closed by assumption. Therefore,
$$\|x_n-x\|'_{V,w}=\|x_n-x\|_{V,w}+\|Tx_n-Tx\|_{W,w'}=\|x_n-x\|_{V,w}+
\|Tx_n-y\|_{W,w'}\longrightarrow  0$$
as $n\longrightarrow  \infty$, which shows that $(V,\|\cdot\|'_w)$  is also
strongly complete.
Since $\|x\|_{V,w}\leq \|x\|'_{V,w}$, for all $x\in V$, Theorem 4.7 shows that
there exists $w''\in (0,1)$ and $C_{w''}$ such that
$\|Tx\|_{W,w'}\leq
\|x\|_{V,w}+\|Tx\|_{W,w'}=\|x\|'_{V,w''}\leq C_{w''}\|x\|_{V,w}$
for all $x\in V$. Therefore by Definition 3.16 we have
$\|T\|_{(w'',w')}<+\infty$, where
$\|T\|_{(w'',w')}=\sup\{\|Tx\|_{w'}|\|x\|_{w''}\leq 1\}$.
Hence by Theorem 3.18,  $T\in B(V,W)$.
\hfill $\Box$
%%%%%
\paragraph{Theorem 4.11 (Uniform boundedness principle theorem in PN-spaces).}
Let $(V,\nu)$ be a strongly complete PN-space and let $(W,\mu)$ be a PN-space,
and let $\{T_n\}_{n\in\BN}$ be a family of operators $T_n\in B(V,W)$, that
satisfy
$$\sup\{\|T_nx\|_{w'}:T_n\in B(V,W)\}<\infty,$$
for all $x\in V$ and $w'\in (0,1)$, then there exists $w\in (0,1)$  such that
$$\sup_{n\in\BN}\|T_n\|_{(w,w')}<\infty$$
%%%%%
\paragraph{Proof.}
First fix $n\in\BN$ and $w'\in(0,1)$, and  define the set
$$F_n=\{x\in V|\sup\|T_nx\|_{w'}\leq n, \; \forall T_n\in B(V,W)\}.$$
Given any $x\in V$, $\sup_{n\in\BN}\|T_nx\|_{w'}<\infty$, by assumption, hence
there exists an integer $n_x\geq 0$ such that
$\sup_{n\in\BN\cup\{0\}}\|T_nx\|_{w'}\leq n_x$, which means that $x\in
F_{n_{x}}$,
 hence $V=\cup_{n=0}^\infty F_n$. On the other hand $F_n$ is
strongly closed subset of $V$.
Since $V$ is strongly complete by Corollary 4.4 (b), there exists an integer
$n_0\geq 0$ such that  $\text{int}
 F_{n_0}\neq \phi$. Hence there exists  $x_0\in
F_{n_0}$ and $r>0$ such that $\overline{B_{V,w}(x_0;r)}\subset F_{n_0}$, where
$$B_{V,w}(x_0;r)=\{x\in V|\|x-x_0\|_w<r\}, \text{for some}~ w\in (0,1).$$
 By definition of $F_{n_0}$, and  Theorem 3.18 for all $z\in
\overline{B_{V,w}(x_0;r)}$ and for all $n\in \BN\cup\{0\}$
% there exists $w'\in (0,1)$:
 $\|T_nz\|_{w'}\leq n_0$
% for all $w'\in (0,1).$
Since any nonzero $x\in V$ can be written as
$x=\frac{\|x\|_w}{r}(z-x_0)$ with $z=(x_0+\frac{r}{\|x\|_w}x)\in
\overline{B_{V,w}(x_0;r)}.$
By Theorem 3.5 we hae
\begin{align*}
\|T_nx\|_{w'} & \leq \frac{\|x\|_w}{r}(\|T_nz\|_{w'}+\|T_nx_0\|_{w'})\\
& \leq \frac{1}{r}(n_0+\|T_nx_0\|_{w'})\|x\|_w\\
& \leq \frac{1}{r}(n_0+\sup_{n\in\BN\cup\{0\}}\|T_nx_0\|_{w'})\|x\|_w,
\end{align*}
for all  $x\in V$.\\
 Therefore, by definition of $\|T\|_{(w,w')}$;
$\sup_{n\in\BN\cup\{0\}}\|T_n\|_{(w,w')}\leq \frac{1}{r}(n_0+
\sup_{n\in\BN\cup\{0\}}\|T_nx_0\|_{w'})<\infty$, since
$\sup_{n\in\BN\cup\{0\}}\|T_nx_0\|_{w'}<\infty$, by assumption.
\hfill $\Box$
%%%%%%
\paragraph{Theorem 4.12 (Banach-Steinhaus Theorem in PN-spaces).}
Let $(V,\nu)$ be a strongly complete PN-space, let $(W,\mu)$ be a PN-space,
let $\{T_n\}_{n=1}^\infty$ be a family of operators $T_n\in B(V,W)$ such that,
for all $x\in V$,  the sequence $\{T_nx\}_{n=1}^\infty$ strongly converges in
$W$, that is $T_nx\longrightarrow  Tx$, as $n\longrightarrow \infty$, then
$T\in B(V,W).$
%%%%
\paragraph{Proof.}
The linearity of each operators $T_n$, combined with the continuity of the
addition and scalar multiplication shows that the operator $T:V\longrightarrow
W$ defined
by $Tx=\lim_{n\longrightarrow \infty}T_nx$ for each $x\in V$ in linear. The
convergence of
each sequence $\{T_nx\}_{n=1}^\infty$ and Theorem 3.17 implies that
: $\sup_{n\in\BN}\|T_nx\|_{w'}<\infty$, for all $x\in V$ and
$w'\in(0,1).$
By Theorem  4.11 there exist $M>0$ and $w\in(0,1)$ such that
$\|T_n\|_{(w,w')}\leq M$, for all $n\in\BN$. Now by Theorem 3.17 and Theorem
3.5, for all $x\in \overline{B_{V,w}(0;1)}$ and  $n\in\BN$ we have
\begin{align*}
\|Tx\|_{w'} & \leq \|Tx-T_nx\|_{w'}+\|T_nx\|_{w'}\\
& \leq \|Tx-T_nx\|_{w'}+M.
\end{align*}
Hence $\|Tx\|_w'\leq M$, when $n\longrightarrow  \infty$, for all $x\in
\overline{B_{V,w}(0;1)}$ therefore by definition of $\|T\|_{(w,w')}$, we have
$T\in B(V,W)$.
\hfill
$\Box$
%%%%%%%%%%%%%%%%%%%%%%%%%%%%%%%%%%%%%%%%%%%

\end{document}